%% file: simple.tex
\documentclass[twocolumn,final]{siamltex}

\input{preamble}

\begin{document}

\title{Work/Precision Tradeoffs in Continuum Models of Biomolecular Electrostatics}
\author{
  Matthew G. Knepley\thanks{knepley@gmail.com, Department of Computational and Applied Mathematics, Rice University, Houston TX 77005} \and
  Jaydeep P. Bardhan\thanks{j.bardhan@neu.edu, Department of Mechanical and Industrial Engineering, Northeastern University, Boston MA 02115}
}

\maketitle

\begin{abstract}
\input{abstract-short.tex}
\end{abstract}

\input{document}

\section*{Acknowledgments}
JPB acknowledges partial support from the National Institute of
General Medical Sciences (NIGMS) of the National Institute of Health
(NIH) under award number R21GM102642. MGK acknowledges partial support
under U.S. DOE Contract DE-AC02-06CH11357. Both authors thank Bob
Eisenberg for the productive environment at Rush University Medical
Center where this work began, and Jed Brown for the suggestion of work
precision diagrams to evaluate algorithmic performance. 

\bibliographystyle{unsrt}
\bibliography{implicit-review,asme}

\end{document}

%% file: preamble.tex
\usepackage{amsmath}
\usepackage{graphicx}
\usepackage{url}

\newcommand{\mO}{\mathcal{O}}
\newcommand{\vq}{\Vec{q}}
\newcommand{\vr}{\Vec{r}}

\graphicspath{ {./figures/} }

%% file: abstract-short.tex
\textit{The structure and function of biological molecules are strongly influenced by the water and dissolved ions that
surround them.  This aqueous solution (solvent) exerts significant electrostatic forces in response to the biomolecule's
ubiquitous atomic charges and polar chemical groups. In this work, we investigate a simple approach to numerical
calculation of this model using boundary-integral equation (BIE) methods and boundary-element methods (BEM).
Traditional BEM discretizes the protein--solvent boundary into a set of boundary elements, or panels, and the
approximate solution is defined as a weighted combination of basis functions with compact support.  The resulting BEM
matrix then requires integrating singular or near singular functions, which can be slow and challenging to compute.
Here we investigate the accuracy and convergence of a simpler representation, namely modeling the unknown surface charge
distribution as a set of discrete point charges on the surface. We find that at low resolution, point-based BEM is more
accurate than panel-based methods, due to the fact that the protein surface is sampled directly, and can be of
significant value for numerous important calculations that require only moderate accuracy, such as the preliminary
stages of rational drug design and protein engineering}

%% file: document.tex

\section*{INTRODUCTION}

Protein structure and function are determined, in large part, by electrostatic interactions between the protein’s atomic
charges and the surrounding biological solvent, a complex mixture of polar water molecules and dissolved charged
ions~\cite{Sharp1990a}. Biological scientists have traditionally modeled these interactions using macroscopic continuum
models based on the Debye--H\"uckel theory or Poisson--Boltzmann equation~\cite{Sharp1990a,Bardhan12_review}. Mean-field
theories of this type assume that solvent molecules are infinitely small compared to the biomolecule
solute~\cite{Beglov97}, a simplification critical to the theoretical studies using spherical and ellipsoidal models of
protein shapes~\cite{Kirkwood34}. In order to accomodate more complex molecular boundaries, the Boundary Element Method
has been widely used~\cite{Rizzo67,Mukherjee03}.

In this paper, we investigate the tradeoff between work and accuracy for representative variants of the BEM
discretization, and quantify the dependence of this tradeoff on mesh resolution and the geometric complexity of the
molecular boundary. Many important computations require only modest accuracy, such as the preliminary stages of rational
drug design and protein engineering, and therefore could benefit from a method with a favorable ratio of flops to
accuracy in this regime. We generate work-precision diagrams comparing panel and point versions of the boundary element
method using the PETSc libraries.

\section*{MODEL}

We will model the problem of protein solvation using the Polarizable Continuum Model (PCM)~\cite{Miertus81,Cammi95,Cances97,Mennucci10}. Thus we
consider a single protein in an infinite solution. In the exterior region, which we label $I$, the potential $\phi$ obeys
Laplace's equation $\nabla^2 \phi = 0$~\cite{Bardhan12_review}, and the dielectric constant is labeled $\epsilon_I$,
which is often taken to be 80, approximately that of bulk water. For simplicity, we omit the source term used in the
linearized Poisson-Boltzmann formulation, since it will not change the arithmetic complexity of our problem. The protein
interior, labeled $III$, is a low-dielectric medium, with $\epsilon_{III}$ usually between 2 and 8. It contains $Q$
discrete point charges, and the potential satisfies a Poisson equation $\nabla^2 \phi = \sum^Q_k q_k \delta(\vr -
\vr_k)$ where $q_k$ and $\vr_k$ specify the $k$th charge. The boundary $a$ separates the protein region $III$ from the
solvent $I$. The potential is assumed to decay to zero suitably fast as $\vr \to \infty$, and the boundary conditions
are continuity of the potential and displacement field,
\begin{align}
  \phi_I(\vr_a) &= \phi_{III}(\vr_a), \\
  \epsilon_I \frac{\partial\phi_I(\vr_a)}{\partial n} &= \epsilon_{III} \frac{\partial\phi_{III}(\vr_a)}{\partial n}.
\end{align}

If we insert the free space Green's function for the potential into the second boundary condition, we obtain a
Boundary Integral Equation~\cite{Altman09,Bardhan11_Knepley} for the induced surface charge on $a$,
\begin{align}\label{eq:BIE}
  \sigma(\vr) + \hat\epsilon \int_a \frac{\partial}{\partial n(\vr)} \frac{\sigma(\vr') d^2\vr'}{4\pi|\vr-\vr'|} =
  -\hat\epsilon \sum^Q_{k=1} \frac{\partial}{\partial n(\vr)} \frac{q_k}{4\pi|\vr-\vr_k|}.
\end{align}
We can write this more compactly using the adjoint of the double layer operator $K$ and the operator $B$ which
gives the normal electric field due to a unit charge,
\begin{align}
  \left( I + \hat\epsilon K^* \right) \sigma &= B \vq \\
  A \sigma &= B \vq
\end{align}
where $\vq$ is the vector of charge values and we introduce an operator $A$.

The \textit{solvation free energy} is the energy of interaction between the induced surface charge, or solvent
polarization, and the solute charges. In order to calculate solvation energies, we will first find the
\textit{reaction potential} $\psi$, which is the potential due to the induced surface charge $\sigma$. We define an
operator $C$ which gives this mapping,
\begin{align}
  \psi(\vr) = C \sigma(\vr) = \int_a \frac{\sigma(\vr') d^2\vr'}{4\pi\epsilon_I |\vr-\vr'|}.
\end{align}
Now we can write the full expression for the solvation energy $E$ in the simple form
\begin{align}
  E = \vq^T C A^{-1} B \vq = \vq^T L \vq,
\end{align}
where we introduce the reaction potential matrix $L$ which has dimensions $Q \times Q$. We will be interested in the
approximation of $L$ below, since it determines the quantity of interest for our biological problems.

\subsection*{PANEL DISCRETIZATION}

A common approach to the discretization of Eq.~\ref{eq:BIE} is the Galerkin method~\cite{Atkinson97,Pozrikidis}, in
which the residual is not required to vanish pointwise, but rather to be orthogonal to a set of test vectors. In the
limit that the test vectors form a basis for the entire approximation space, we recover the same solution as the
original equation. Thus we introduce test functions $\tau_i(\vr)$, and integrate our equation
\begin{align}\label{eq:panelBEM}
  \int_a d^2\vr\, \tau_i(\vr) \left( I + \hat\epsilon K^* \right) \sigma(\vr) &= \\
  \quad \int_a d^2\vr\, \tau_i(\vr) B \vq \nonumber\\
  \int_a d^2\vr\, \tau_i(\vr) \left( \sigma(\vr) + \hat\epsilon \int_a \frac{\partial}{\partial n(\vr)}
  \frac{\sigma(\vr') d^2\vr'}{4\pi|\vr-\vr'|} \right) &= \\
  \quad -\int_a d^2\vr\, \tau_i(\vr) \hat\epsilon \sum^Q_{k=1} \frac{\partial}{\partial n(\vr)} \frac{q_k}{4\pi|\vr-\vr_k|}. \nonumber
\end{align}

We must decide how to discretize the integrals appearing in Eq.~\ref{eq:panelBEM}, using some sort of quadrature. Below
we employ triangular panels and constant basis functions in each panel. For quadrature, we choose the
\textit{qualocation} strategy~\cite{Tausch01,Bardhan09_disc}, which uses a single point centroid quadrature 
for the inner integration over $\vr'$, and either a high-order quadrature or analytic result for the outer integration
over $\vr$,
\begin{align}\label{eq:panelBEMQual}
  \sum_m w_m \tau_i(\vr_m) \left( \sigma(\vr_m) + \hat\epsilon \frac{\partial}{\partial n(\vr_m)}
  \frac{\sigma(\vr_j) a_j}{4\pi|\vr_m-\vr_j|} \right) &= \\
  \quad -\sum_m w_m \tau_i(\vr_m) \hat\epsilon \sum^Q_{k=1} \frac{\partial}{\partial n(\vr_m)} \frac{q_k}{4\pi|\vr_m-\vr_k|}. \nonumber
\end{align}
In the results section, we use analytic integrals, either the Hess-Smith method, or the Newman method when
vector to the evaluation point is too close to the panel normal~\cite{Hess62,Newman86}.

In order to compute one entry in our matrices $A$ or $B$, we need to integrate over a panel and perform an evaluation at
our source point. In the following, we will only consider the cost of evaluating $A$ since $B$ and $C$ are much less
expensive due to the fact the $Q << N_c$, where $N_c$ is the number of panels on the surface. Once we have our matrices
$A$, $B$, and $C$, the cost to compute the $L$ matrix is dominated by the inversion of $A$, which executes $\mO(N_c^3)$
flops. If we let the relative cost of panel integration against point evaluation be $N_p$, then the dominant terms in
our complexity model are
\begin{align}
  \mO(N_p N_c^2) + \mO(N_c^3).
\end{align}
Clearly, as the size of the mesh grows, the time will be completely dominated by the second term. However, this is
usually alleviated by using an iterative method for the inversion, which can reduce the cost to
\begin{align}\label{eq:panelCoplexity}
  \mO(N_p N_c^2) + \mO(N_c^2)
\end{align}
since the number of iterates is bounded by a constant independent of the size of the system because the condition number
of $A$ is bounded~\cite{Atkinson97}. It is also true that the exponent of $N_c$ can be reduced from 2 to 1 using methods
that exploit the decay in the interaction with distance, such as the Fast Multipole Method~\cite{Greengard88}, modified
FFT methods~\cite{Altman2006}, or skeletonization~\cite{Ho12_Greengard}. This will not be our main concern, since we are evaluating cost $N_p$ incurred by panel integration against the extra accuracy gained.

\subsection*{POINT DISCRETIZATION}

As a lower-cost alternative to the panel discretization, we consider a collocation method which enforces the continuum
equation at a given set of points $\{\vr_i\}$, and discretize the integral operator using the Nystr\"om method (quadrature),
\begin{align}
  \left( I + \hat\epsilon K^* \right) \sigma(\vr_i) &= B \vq \\
  \sigma(\vr_i) + \hat\epsilon \int_a \frac{\partial}{\partial n(\vr_i)} \frac{\sigma(\vr') d^2\vr'}{4\pi|\vr_i-\vr'|} &=
  -\hat\epsilon \sum^Q_{k=1} \frac{\partial}{\partial n(\vr_i)} \frac{q_k}{4\pi|\vr_i-\vr_k|} \\
  \sigma(\vr_i) + \hat\epsilon \sum_j \frac{\partial}{\partial n(\vr_i)} \frac{\sigma(\vr_j) w_j}{4\pi|\vr_i-\vr_j|} &=
  -\hat\epsilon \sum^Q_{k=1} \frac{\partial}{\partial n(\vr_i)} \frac{q_k}{4\pi|\vr_i-\vr_k|}
\end{align}
and the contribution from the $i = j$ term is 0. A similar collocation method is used to evaluate the $C$ operator. For
this study, we have chosen to use the vertices of panel triangulations as the points, and the point weights are
determined by giving each vertex 1/3 of the panel area for each panel incident to it.

The construction of $A$, $B$, $C$ and $L$ for the point discretization has a complexity model very similar to
Eq.~\ref{eq:panelCoplexity},
\begin{align}
  \mO(N_c^2) + \mO(N_c^2).
\end{align}

\section*{RESULTS}

We demonstrate the correctness of our implementation using a simple problem with known exact solution. Then we show that
our point BEM variants can be more efficient in a certain accuracy range, using realistic test cases for amino acid
residues.  For all test cases, we take the dielectric constant of the solvent $\epsilon_I = 80$, and the dielectric constant of the
solute $\epsilon_{III} = 4$. Initial MATLAB implementations of the panel BEM method and point BEM
method are available in online repositories, along with a higher-performance implementation using the PETSc
libraries~\cite{petsc-web-page,petsc-user-ref} was used to run all the tests below, and is available
online~\cite{PointBEMPETScRepo}. All timings are from an Apple Airbook with a 1.4 GHz Core i5 processor and 8GB
of 1600 MHZ DDR3 RAM.

\subsection*{SPHERE}

As an initial test, we compare the error in the methods described above for the problem of a low-dielectric sphere in a
high-dielectric medium, which has a classical exact solution~\cite{Kirkwood34,Bardhan11_Knepley}. Specifically, we take a sphere of radius
$R$ with dielectric constant $\epsilon_{III}$, embedded in a medium of dielectric constant $\epsilon_I$, and fill it
with random charges. The charges are randomly selected from vertices on a grid of spacing $h$, where selected vertices
must be at least distance $h$ from the spherical surface.

The reaction potential $\psi$, which we evaluate at each charge location $x_j$, is given by
\begin{equation}\label{eq:ReactionPotentialExpansion}
  \psi = \sum^\infty_{n=0} \sum^n_{m=-n} B_{nm} r^n_j P_n^m(\cos\theta_j) e^{i m \phi_j}.
\end{equation}
where the reaction potential expansion coefficients are given by
\begin{equation}\label{eq:Bnm}
  B_{nm} = \frac{1}{\epsilon_1 b^{2n+1} } \frac{(\epsilon_1 - \epsilon_2) (n+1)}{\epsilon_1 n + \epsilon_2 (n+1)} E_{nm},
\end{equation}
and the charge distribution expansion coefficients are
\begin{equation}
  E_{nm} = \sum^Q_{k=1} q_k r^n_k \frac{(n-|m|)!}{(n+|m|)!} P_n^m(\cos\theta_k) e^{-i m \phi_k}.
\end{equation}

\begin{figure}[t]
\begin{center}
\includegraphics[width=0.5\textwidth]{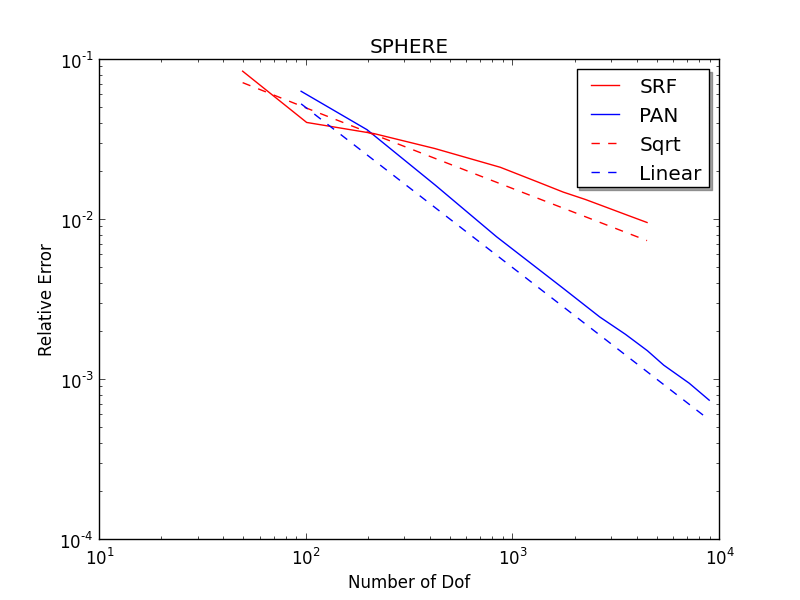}
\end{center}
\caption{WE SHOW MESH CONVERGENCE FOR BEM METHODS APPLIED TO THE PROBLEM OF RANDOM CHARGES IN A SPHERE, WITH $R = 6.0$, $h = 1.0$, $\epsilon_I
  = 80$, $\epsilon_{III} = 4$, AND $Q = 10$ CHARGES. THE ANALYTIC SOLUTION IS CALCULATED TO MULTIPOLE ORDER 25. CLEARLY,
THE PANEL METHOD (PAN) CONVERGES LINEARLY, AND THE POINT METHOD SRF HAS ORDER 1/2.}
\label{fig:meshConvSphere} 
\end{figure}

In order to demonstrate correctness, we show convergence of our  solution to the analytic reference
solution as the element size is reduced, and call this \textit{mesh convergence}. In Fig.~\ref{fig:meshConvSphere}, we
demonstrate mesh convergence of both BEM variants for the sphere problem. As we expect, the panel method converges
with order 1, and the point method with order 1/2~\cite{Atkinson97}.

\begin{figure}[t]
\begin{center}
\includegraphics[width=0.5\textwidth]{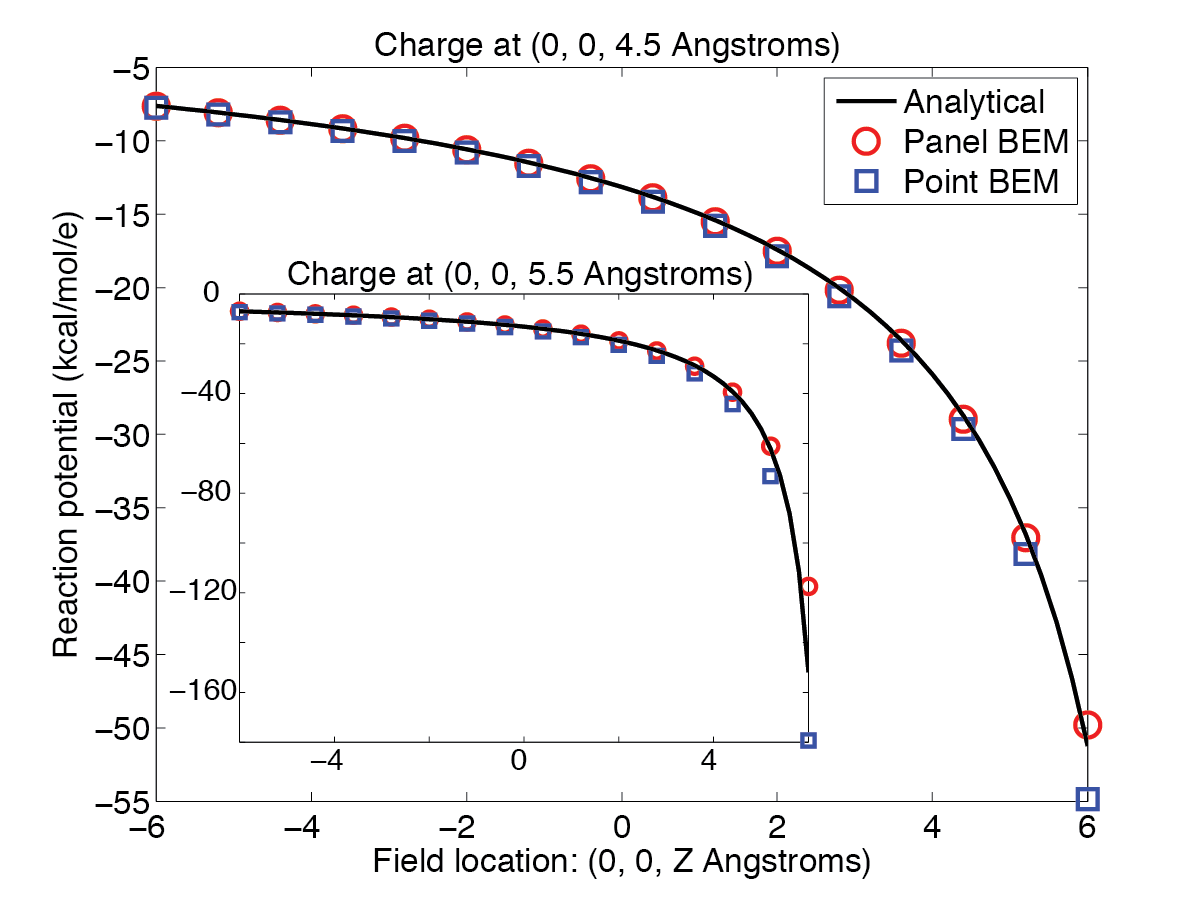}
\end{center}
\caption{WE SHOW THE RESULTING POTENTIALS FOR SINGLE CHARGES IN A SPHERE, WITH $R = 6.0$, $h = 1.0$, $\epsilon_I
  = 80$, $\epsilon_{III} = 4$. THE CHARGE IS AT (0, 0, 4.5 ANGSTROMS) FOR THE MAIN PLOT AND AT (0, 0, 5.5 ANGSTROM) FOR THE INSET.   POINT BEM ACCURACY IS ACCEPTABLE FOR THE CHARGE AT 4.5 ANGSTROMS BUT NOT FOR THE ONE AT 5.5 ANGSTROMS. }
\label{fig:lineReactionPotentials} 
\end{figure}

To show a practical case, we note that although  accurate representation of the boundary is important globally, protein charges are at atom centers and therefore usually more than about 1 Angstrom from the boundary.  Fig.~\ref{fig:lineReactionPotentials} shows the potentials for single charges in the 6-Angstrom sphere: for the main plot, the charge is at $(0, 0, 4.5~\mathrm{\AA})$ and $(0, 0, 5.5~\mathrm{\AA})$ for the inset.  We plot potentials along the $Z$ axis.  For the $Z=4.5$ charge, the accuracy is acceptable throughout  (less than 1~kcal/mol/$e$), except between the charge and the boundary.  On the other hand, the error is large for the less realistic case (note that the scales are different).

\begin{figure}[t]
\begin{center}
\includegraphics[width=0.5\textwidth]{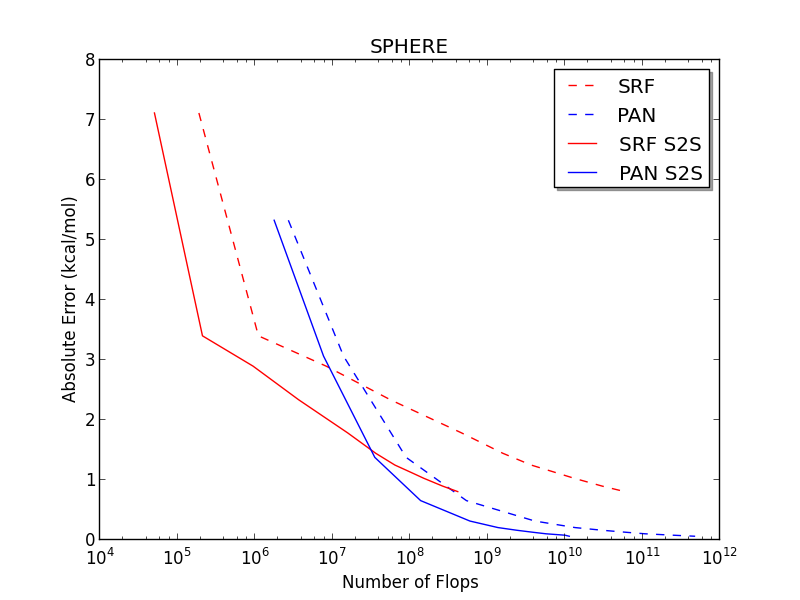}
\end{center}
\caption{WE SHOW A WORK-PRECISION DIAGRAM FOR BEM METHODS APPLIED TO THE SAME PROBLEM AS
  FIG.~\ref{fig:meshConvSphere}. FOR THIS CASE, THE PANEL METHOD IS SUPERIOR FOR ERROR BELOW 3 KCAL/MOL, ALTHOUGH WHEN
  ONLY CONSIDERING THE SURFACE-TO-SURFACE (S2S) OPERATOR $A$, THE POINT METHOD IS SUPERIOR DOWN TO 1.5 KCAL/MOL.}
\label{fig:workPrecSphere} 
\end{figure}

In order to look more closely at the computation itself, we will employ a work-precision diagram~\cite{Hairer2009}. This plots the work done
against the accuracy that was achieved by the computation. We define the work as the number of floating point operations
(flops) used to compute a) the reaction potential matrix $L$, and b) just the entries of $A$. We expect the second metric to map
well to the case where scalable, or fast, algorithms are used for construction since the local work done for those
algorithms is typically greater than half of the total, and uses exactly these routines, whereas the large number of
flops in the direct inversion of $A$ dominates the total for the simpler dense case. When we examine
Fig.~\ref{fig:workPrecSphere}, we see that to obtain errors less than 3 kcal/mol it is more advantageous to use the
panel method. If we only consider the formation of the surface-to-surface operator $A$, this accuracy threshold
decreases to 1.5 kcal/mol.

\subsection*{RESIDUES}

To compare the discretizations on realistic, atomically detailed molecular geometries, we next look at the reaction field of two $\alpha$-amino acids, aspartic acid (ASP) and arginine (ARG).  Amino acids are the building blocks of proteins, and because these and other charged amino acids play key roles in biomolecular electrostatics, accurate calculations are essential. Meshes were created for
the isolated amino acid structures using the meshmaker program from the FFTSVD package~\cite{Altman06}, the MSMS program~\cite{Sanner96},
PDB files from the Protein Data Bank~\cite{Berman00}, and atomic radii parameterized for the CHARMM force field~\cite{Brooks83,Nina97}. Using
these tools, we can create a range of panel/point densities.

\begin{figure}[t]
\begin{center}
\includegraphics[width=0.5\textwidth]{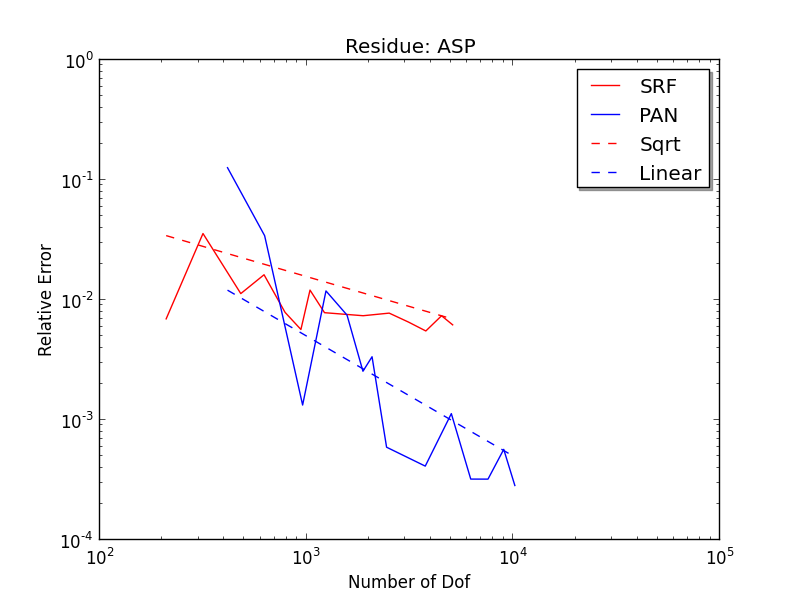}
\end{center}
\caption{WE SHOW MESH CONVERGENCE FOR BEM METHODS APPLIED TO THE ASPARTIC ACID MOLECULE, WITH $\epsilon_I = 80$ AND
  $\epsilon_{III} = 4$. WE AGAIN SEE THAT THE PANEL METHOD (PAN) CONVERGES LINEARLY, AND THE POINT METHOD (SRF) HAS
  ORDER 1/2, ALTHOUGH THE CONVERGENCE IS INITIALLY NOISY.}
\label{fig:meshConvASP} 
\end{figure}

\begin{figure}[t]
\begin{center}
\includegraphics[width=0.5\textwidth]{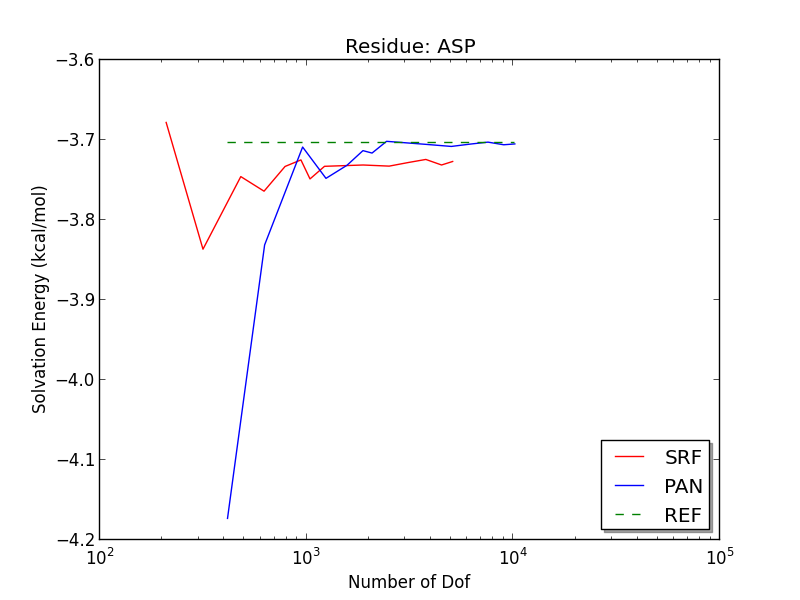}
\end{center}
\caption{WE SHOW SOLVATION ENERGY CONVERGENCE FOR BEM METHODS APPLIED TO THE ASPARTIC ACID MOLECULE ($\epsilon_I = 80$ and $\epsilon_{III} = 4$). NOTE THAT INITIALLY THE POINT METHOD CONVERGES FASTER, BUT STAGNATES SOON AFTER.}
\label{fig:meshConvASPEnergy} 
\end{figure}

\begin{figure}[t]
\begin{center}
\includegraphics[width=0.5\textwidth]{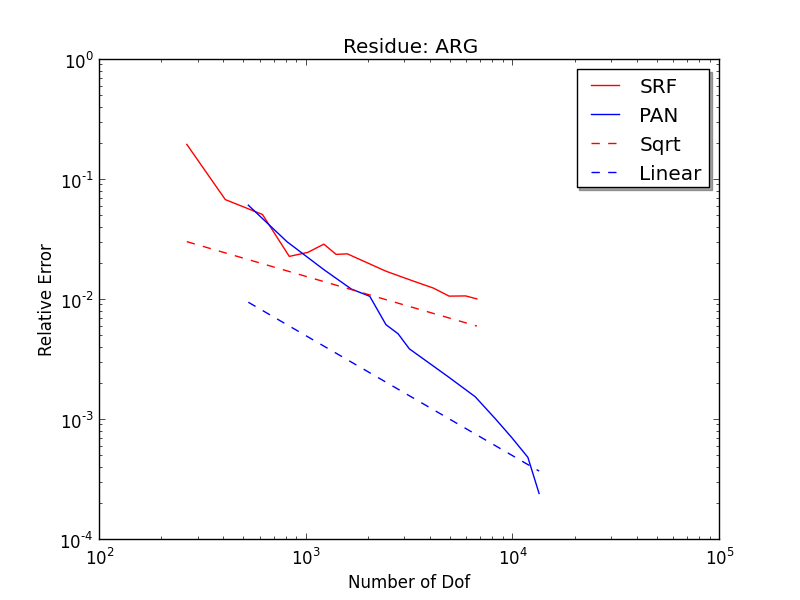}
\end{center}
\caption{WE SHOW MESH CONVERGENCE FOR BEM METHODS APPLIED TO THE ARGININE MOLECULE ($\epsilon_I = 80$ AND
  $\epsilon_{III} = 4$).  AGAIN THE PANEL METHOD (PAN) CONVERGES LINEARLY, AND THE POINT METHOD (SRF) HAS
  ORDER 1/2, ALTHOUGH  CONVERGENCE IS INITIALLY NOISY.}
\label{fig:meshConvARG}
\end{figure}

In Fig.~\ref{fig:meshConvASP}, we show mesh convergence for ASP. It is initially noisy, but we clearly see the
asymptotic convergence rates for larger meshes. We have defined the reference energy to be the result of Richardson
extrapolation~\cite{Richarson1911} using the last two energies from the panel method. Thus, in all residue figures the
last panel method energy appears unrealistically accurate. The absolute solvation energies are shown in
Fig.~\ref{fig:meshConvASPEnergy}, and it is clear that the order 1/2 convergence of the point method looks like
complete stagnation on the scale of our tests. The same convergence behavior is seen in Fig.~\ref{fig:meshConvARG} for ARG.

\begin{figure}[t]
\begin{center}
\includegraphics[width=0.5\textwidth]{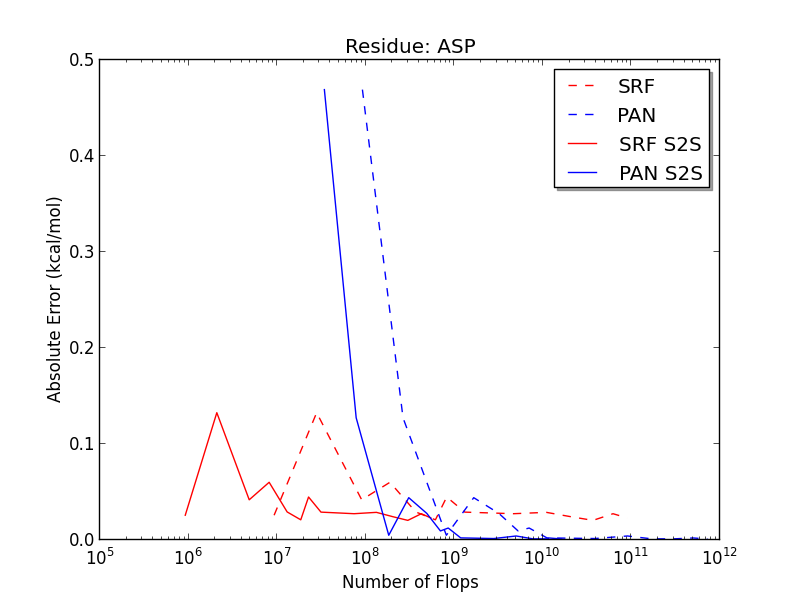}
\end{center}
\caption{WE SHOW A WORK-PRECISION DIAGRAM FOR BEM METHODS APPLIED TO THE SAME PROBLEM AS
  FIG.~\ref{fig:meshConvASP}. FOR THIS CASE, THE PANEL METHOD IS SUPERIOR FOR ERROR BELOW 0.25 KCAL/MOL.}
\label{fig:workPrecASP}
\end{figure}

\begin{figure}[t]
\begin{center}
\includegraphics[width=0.5\textwidth]{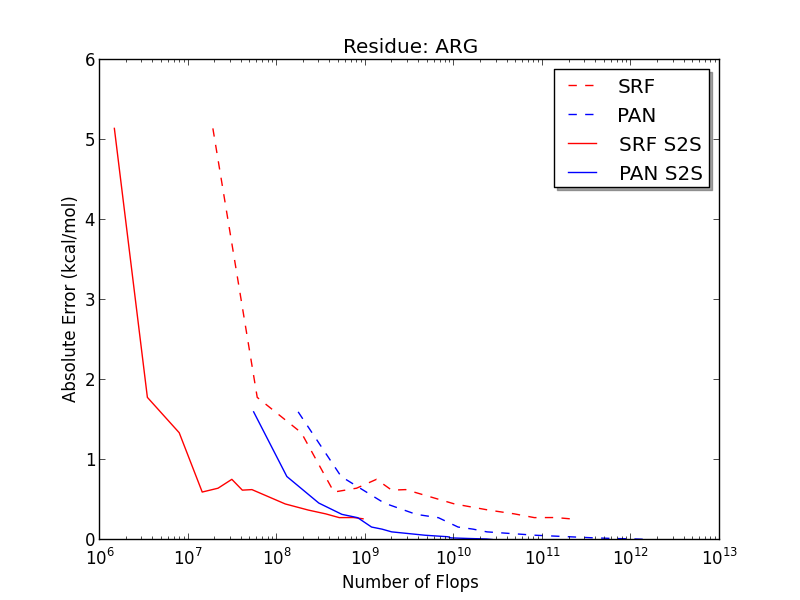}
\end{center}
\caption{WE SHOW A WORK-PRECISION DIAGRAM FOR BEM METHODS APPLIED TO THE SAME PROBLEM AS
  FIG.~\ref{fig:meshConvARG}. FOR THIS CASE, THE PANEL METHOD IS SUPERIOR FOR ERROR BELOW 0.75 KCAL/MOL, ALTHOUGH WHEN
  ONLY CONSIDERING THE SURFACE-TO-SURFACE (S2S) OPERATOR $A$, THE POINT METHOD IS SUPERIOR DOWN TO 0.25 KCAL/MOL.}
\label{fig:workPrecARG} 
\end{figure}

From the work-precision diagrams in Fig.~\ref{fig:workPrecASP} and Fig.~\ref{fig:workPrecARG}, we see that the accuracy
threshold for the point method to be preferred has moved considerably lower. In fact, when comparing construction of the
surface-to-surface operator $A$, we see that the point method has lower costs up to errors of 0.25 kcal/mol, well within
the range of many biological studies.

\section*{DISCUSSION}

We have compared the work-accuracy tradeoff in panel and point BEM methods. Very clear accuracy thresholds emerge below
which point BEM is preferable to the panel method. Point BEM can be used productively for scenarios that tolerate low
resolution, such as intermediate calculations in design iterations, structure optimization, and optimization of binding
affinity. The range of applicable problems could be expanded by optimizing the placement and weighting of
points. Sophisticated tools, such as MSP~\cite{MSP}, could be used to improve our point discretizations, and is the
subject of future research.

How does our outlook change with the introduction of fast methods for the application of the integral operators? The
characteristics of the local work remain identical. On modern architectures, an optimal workload roughly balances long
range, or tree work, with local direct calculations~\cite{GreengardGropp1990}. Thus, the speedup we show above of point
methods over panel methods should remain roughly intact. We plan to validate this claim using the PyGBe fast
solver~\cite{Cooper13}.